\documentclass[12pt,titlepage,a4paper]{article}
\usepackage{amssymb,latexsym,amsfonts,amsmath,amsthm}

{\makeatletter \@addtoreset{equation}{section}}

\newtheorem{theorem}{Theorem}[section]
\newtheorem{lemma}[theorem]{Lemma}

\newcommand{\blackbox}{\hfill \rule{2mm}{2mm}}
\newcommand{\diag}{\mathop{\mathrm{diag}}\nolimits}
\newcommand{\Ker}{\mathop{\mathrm{Ker}}\nolimits}
\newcommand{\IIm}{\mathop{\mathrm{Im}}\nolimits}

\begin{document}
\title{Homomorphisms of modules associated with polynomial matrices with
infinite elementary divisors}

\author{Pudji Astuti \\
  Departement of Mathematics \\
and Natural Sciences\\
Institut Teknologi Bandung\\
Bandung 40132\\
Indonesia
\thanks{The research of the first author was supported by
 Deutscher Akademischer Austauschdienst under Award No.~A/98/25636.} 
\and
Harald K. Wimmer\\
Mathematisches Institut\\
Universit\"at W\"urzburg\\
D-97074 W\"urzburg\\
Germany}

\date{}

\maketitle

\begin{abstract}
If the inverse of a nonsingular polynomial matrix $L$
 has a polynomial part
then one can associate with $L$ a module over the ring of proper rational
functions, which is related to  the structure of L at infinity.
 In this paper we characterize homomorphisms of such modules.

\vspace{2cm}

\vspace{.5cm}

\end{abstract}

\section{Introduction}

According to Rosenbrock \cite{Rbr} a transfer matrix $G \in
K^{m \times p}$ of rational functions over a field $K$ admits
a \emph{generalized state space realization}
\[
 G(s) = \begin{pmatrix} C_1 & C_2 \end{pmatrix}
      \begin{pmatrix}
                        sI - A_1 & 0 \\
                 0 & s N_2 - I
       \end{pmatrix}
                       \begin{pmatrix} B_1 \\ B_2
          \end{pmatrix}
\]
       such that
\begin{equation} \label{e.prr}
     G_1(s) = C_1 (sI - A_1) ^{-1} B_1
\end{equation}
is the strictly proper part 
and
\begin{equation} \label{e.pol}
 G_2(s) = C_2 (s N_2 - I) ^{-1} B_2,
\end{equation}
where $N_2$ is nilpotent, is the
 polyonomial part of $G$.
It is well known that the realizations \eqref{e.prr} and
\eqref{e.pol} can be constructed by module theoretic approaches.
In the case of \eqref{e.prr} a construction is due to
Fuhrmann \cite{Fu2}. For a realization
theory of anticausal input
output maps we refer to Conte and Perdon \cite{CP1}.
To describe the polynomial models that serve as state spaces for
\eqref{e.prr} and  \eqref{e.pol} we use the following notation.
A rational function $f \in K(s)$   is called \emph{proper}
or \emph{causal} (resp. \emph{strictly proper} or
\emph{strictly causal}) if
$f = 0$ or if $f \neq 0$
and $f  = p/q, \, p,q \in K[s], q \neq 0$,
and ${\rm{deg}}\, p \leq \textrm{deg}\, q$
(resp.  ${\rm{deg}}\, p < \textrm{deg}\, q$).
Let $K_{\infty}(s)$ denote the ring of proper rational
functions over $K$. Then
\begin{equation} \label{e.oplus}
         K(s) = K[s] \oplus  s^{-1} K_{\infty}(s).
\end{equation}
To \eqref{e.oplus} correspond projection operators
\[
  \pi _- :   K(s) \to  s^{-1} K_{\infty}(s)
\]
and
\[
     \pi _+  = (I - \pi _- ):   K(s) \to  K[s].
\]
Put
\begin{equation} \label{e.f0}
           (f)_0 = (\pi _+ f)(0), \,  f \in K(s).
\end{equation}
The decomposition \eqref{e.oplus}, the projections
$ \pi _-$ and $ \pi _+$, and  definition
\eqref{e.f0} extend naturally from $K(s)$
 to  $K^n(s)$
and $K^{m \times p}(s)$.

Let $G \in K^{m \times p}(s)$ have a realization
\begin{equation} \label{e.rlz1}
        G = W_1 + P_1 D_1 ^{-1} Q_1
\end{equation}
where  $W_1, P_1,  Q_1, D_1$ are polynomial
matrices, with $D_1$ of size $n_1 \times n_1 $.
In Fuhrmann's theory \cite{F1} a state space for a realization
\eqref{e.prr} of $\pi _{-} G$ is provided by
\[
   V _{D_1} =  K^n_{1}[s]/D_1 K^n_{1}[s].
\]
Obviously $ V _{D_1}$ is a $K[s]$-module and therefore
also a vector space over  $K$.
The counterpart of \eqref{e.rlz1}
is a realization
\begin{equation} \label{e.rlz2}
        G = W_2 + P_2 D_2 ^{-1} Q_2,
\end{equation}
where
$ P_2 $ and $Q_2$ are proper rational matrices,
$ W_2 $    is strictly proper rational and
$ D_2 $ is a polynomial matrix,
$ D_2 \in  K^{n_{2} \times n_{2}} $.
Define
\begin{equation} \label{e.ud2}
  U ^{D_2} =
 K^{n_{2}}_{\infty}(s)/\bigl( K^{n_{2}}_{\infty}(s)
    \cap D_2 s^{-1}  K^{n_{2}}_{\infty}(s)\bigr).
\end{equation}
Then $U ^{D_2}$ is a $K_{\infty}(s)$-module  and at the
same time a $K$-vector space. At the end of this section
we shall indicate why  $U ^{D_2}$ can be taken as a state
space of a  realization \eqref{e.pol} of $\pi _{+} G$.
 Let us mention that the finite
and infinite \emph{pole modules} (see \cite{WSCP}) of $G(s)$ are given
by $V_{D_1}$ and $ U ^{D_2}$, if \eqref{e.rlz1} is an irreducible
realization and \eqref{e.rlz2} satisfies coprimeness conditions
of the form \eqref{e.left}.

We note that a nonsingular polynomial matrix $L \in
K^{n \times n}[s]$ gives rise to two types of
 modules, namely the
 $K[s]$-module
\[
     V_L = K^n[s]/L K^n[s]
\]
and the  $K_{\infty}(s)$-module
\begin{equation} \label{e.UL}
   U^{L} =
 K^n_{\infty}(s)/\bigl( K^n_{\infty}(s)
\cap L s^{-1} K^n_{\infty}(s)\bigr).
\end{equation}
Beside realizations there is a wide range of issues such as similarity of state space
models, system equivalence or simulation of restricted input output
maps  which involve two polynomial matrices
$L$ and $L_1$ and  homomorphisms from $V_L$ to $V_{L_1}$ and
 from $U^L$ to $U^{L_1}$.
The  $K[s]$-module
homomorphisms from  $V_L$ to $V_{L_1}$  are well understood.
According to
Fuhrmann \cite{F1} their description is based on
intertwining relations between $L$ and $L_1$.
 In this note we study  $K_{\infty}(s)$-module
homomorphisms from  $U^L$ to  $U^{L_1}$. Our characterizations will be in
correspondance with Fuhrmann's results in Ref. 
  \cite{Fu2,F1}. 
  Comparing the definitions of $V_L$ and
$U^L$ we observe that $L K^n[s]$ is a submodule of $K^n[s]$
whereas in general $Ls^{-1}K^n_{\infty}(s)$ is not contained in
$K^n_{\infty}(s)$. Hence  it is not surprising that
$U^L$ is less easy to handle than $V_L$ and that in our study
 technical obstacles have to be removed  which do not appear
in the case of the module $V_L$.

To obtain a concrete representation of $U^L$   we define
a map
\[
   \rho  ^L : K^n_{\infty}(s)  \to K^n[s]
\]
      by
\[
         \rho  ^L x = L \pi _+ L^{-1} x, \,\,  x \in K^n_{\infty}(s).
\]
Put $\bar{x} = \rho  ^L x$. For $q \in K_{\infty}(s)$ and
 $\bar{x} \in \IIm  \rho ^L$ we set
    $   q \cdot \bar{x} = \overline{qx}. $
This product is well defined since
\[
       \Ker \rho ^L =
\bigl( K^n_{\infty}(s) \cap s^{-1} L K^n_{\infty}(s) \bigr).
\]
Therefore $\IIm \rho ^L$ is a $ K_{\infty}(s)$-module, isomorphic
to the
quotient module $U^L$ in \eqref{e.UL}. From now on we identify both
modules such that
\[
   U^L = \IIm \rho ^L =  L \pi _+ L^{-1}K^n_{\infty}(s).
\]
Clearly, $U^L = 0$ \, if \,$sL^{-1}$ \, is  proper rational.
A shift operator $S_{-}(L)$ on $U^L$ is given by
\[
    S_{-}(L) \bar{x} =   s^{-1}  \cdot \bar{x}, \,\,
          \bar{x} \in U^L .
\]
Clearly, $S_{-}(L)$ is a nilpotent endomorphism of  $U^L$.

Let us now give a concrete example for the use
of $ K_{\infty}(s)$-module $U^L$.Based on the representation \eqref{e.rlz2}
of $ G$  we derive a  realization of $\pi _{+} G$
having  $U^{D_2}$ as its state space. We adapt a construction of
\cite{Fu3}. Assume $\pi _{+} G(s) = \sum ^{t}_{\nu = 0}
G _{\nu}  s^{\nu}$. Define the map
$B_2 : K^p \to U^{D_2}$ by
\[
          B_2 \xi = \rho ^{D_2}\, Q_2 \xi, \,\,  \xi \in K^p.
\]
Put $ N_2 = S_{-} (D_2)$ and
define $C_2 : U^{D_2} \to K^m$ by
\[
        C_2  \bar{x} = - \left( P_2 D_2 ^{-1} \bar{x} \right)_0,
       \,\,  \bar{x} \in  U^{D_2}.
\]
Then a straightforward computation yields
\[  G _{\nu} = -  C_2  N_2 ^{\nu} B_2, \,\, \nu = 0,1, \dots , t,
\]
such that
\[
     \sum ^{t}_{\nu = 0} G _{\nu}  s^{\nu}  =
        C_2 (s N_2 - I) ^{-1}  B_2.
\]

\section{Basic facts of  
the module {\boldmath{$U^L$}}}

For a nonzero proper rational function
$f = p/q , \,\, p,q \in K[s]$,
let  a degree function be defined by
$\delta (p/q) = {\rm{deg}}q - \textrm{deg}p$. It is well known that
 $\bigl( K_{\infty}(s), \delta \bigr)$
is a euclidean domain.
The units $K^{*} _{\infty}(s)$  are the proper rational functions
$f$ with $\delta f = 0$. The ideal $(s^{-1})$ is the unique maximal
ideal of $K_{\infty}(s)$. Let us call
a matrix $P \in K^{n \times n} _{\infty}(s)$ \emph{bicausal}
if $\det P   \in K^{*}_{\infty}(s)$, i.e. if $P$ is invertible in
$ K^{n \times n} _{\infty}(s)$.
If $W \in K^{m \times r}(s)$ has rank $n$ then there exist bicausal matrices
$P$ and $Q$ such that
\[
     W = P \left( \begin{matrix}
                                  \Sigma  & 0 \\
                              0 & 0
           \end{matrix} \right) Q
\]
with
\begin{multline} \label{e.SMM}
  \Sigma = \diag (s^{-\alpha_1}, \dots , s^{-\alpha_t}, s^{\beta_{t+1}},
\dots ,  s^{\beta_n}),\\
     -\alpha_1 \leq \dots \leq -\alpha_t < 0
 \leq \beta_{t+1} \leq
\cdots \beta_n.
\end{multline}
The integers $ -\alpha_1, \dots , \beta_n$ are uniquely determined by $W$.
In particular, if  $L \in K^{n \times n}[s]$ is nonsingular then
\begin{equation} \label{e.infS}
      s^{-1} L = P \Sigma Q
\end{equation}
for some $P,Q \in K^{n \times n}_{\infty}(s)^{*}$ and $\Sigma$ as in
\eqref{e.SMM}.  In the case of a linear pencil $L(s) = A_0 - A_1s$\, the
polynomials
   $ s^{\alpha _1}, \dots , s^{\alpha _t} $
 are the elementary divisors of
$A_0s - A_1$ belonging to the characteristic root $0$.
According to \cite{Wi1} the matrix $\Sigma$
in \eqref{e.infS} and \eqref{e.SMM} provides information on the structure
 of $U^L$.
We have
\[
    U^L \cong \oplus
      \bigl\{  K_{\infty}(s) /  s^{-\alpha _j} K_{\infty}(s), \,
    j = 1, \dots , t \bigr\}
\]
such that $U^L$ is a finitely generated torsion module over  $K_{\infty}(s)$
with elementary divisors
\begin{equation} \label{e.ueninf}
          s^{- \alpha _1}, \dots ,  s^{- \alpha _t}.
\end{equation}
We call \eqref{e.ueninf} the \emph{infinite elementary divisors}
of $L$. Then $ s^{ \alpha _1}, \dots ,  s^{ \alpha _t}$ are the
elementary divisors of the shift $S_{-}(L)$, and $\mathrm{dim}_K\,U^L
= \alpha _1 + \dots +  \alpha _t$.
To describe a dual pairing  \cite{Wi2}
 between
the $K$-linear spaces $U^{L^{T} }$ and $U^L$
we note   that
\begin{equation} \label{e.dual}
   \langle \bar{y}, \bar{x} \rangle = (y^T L^{-1} x)_0, \, \,
   \bar{y} \in U^{L^{T} },\, \bar{x} \in U^L,
\end{equation}
is  a well defined nondegenerate bilinear form on
$U^{L^{T} } \times U^L$.

\section{Homomorphisms}

Our main result is Theorem \ref{Th.main} below. Its proof will be based on
the subsequent two lemmas. In the following
$L \in  K^{n \times n}_{\infty}(s)$ and
$L_1 \in K^{n_1 \times n_1}_{\infty}(s)$ will be fixed nonsingular polynomial
matrices.

\begin{lemma} \label{La.toL1}
A map
\begin{equation} \label{e.Phi}
          \Phi : K^{n}_{\infty}(s) \to U^{L_1}
\end{equation}
is a $ K_{\infty}(s)$-module homomorphism if and only if there exists
a matrix
             $\Theta \in  K^{n_1 \times n}_{\infty}(s)$
such that
\begin{equation} \label{e.Phx}
          \Phi x = \rho ^{L_1} (\Theta x), \,\,  x \in  K^{n}_{\infty}(s).
\end{equation}
\end{lemma}

\emph{Proof.} Let $e_1, \dots , e_n$ be the standard basis of $K^n$. Assume
that $\Phi$ in \eqref{e.Phi} is a  $ K_{\infty}(s)$-module homomorphism.
Then $\Phi e_i = \rho ^{L_1} \theta _i$ for some
$\theta _i \in K^{n_1}_{\infty}(s)$ and \eqref{e.Phx} holds with
$\Theta = ( \theta _1, \dots \theta _n)$. The converse is obvious. \blackbox

\vspace{2.5mm}
Condition \eqref{e.kerincl}  below
together with a somewhat technical equivalent condition   will be crucial.

\begin{lemma} \label{La.TT1}
We have
\begin{equation} \label{e.kerincl}
   \Theta \Ker \rho ^L \subseteq \Ker \rho ^{L_1 }.
\end{equation}
with $\Theta \in K^{n_1 \times n}_{\infty}(s)$ if and only if
there exist a matrix
$\Theta _1 \in K^{n_1 \times n}_{\infty}(s)$  and a matrix $\Psi$
satisfying
\begin{equation} \label{e.Psi}
\Psi \in s^{-1}K^{n_1 \times n}_{\infty}(s) \,\,\, and \,\,\,
 L_1 \Psi \in K^{n_1 \times n}_{\infty}(s)
\end{equation}
such that
\begin{equation} \label{e.L1Psi}
    (\Theta +  L_1 \Psi) L = L_1 \Theta _1.
\end{equation}
\end{lemma}

\emph{Proof.}  It is evident that \eqref{e.L1Psi} implies \eqref{e.kerincl}.
To prove the converse implication we note that \eqref{e.kerincl}
is equivalent to \,$
   \Theta \Ker \rho ^L \subseteq s^{-1}L_1K^{n_1}_{\infty}(s).$
If \linebreak
 $s^{-1}L$ is factorized as in \eqref{e.infS},
\begin{multline} \label{e.AB}
        s^{-1} L = P \Sigma Q, \,\, \Sigma = \diag (A, \, B), \, \\
A =\diag (s^{-\alpha _1}, \dots , s^{-\alpha _t}), \,
B =\diag (s^{\beta _{t+1} }, \dots ,  s^{\beta _n})
\end{multline}
then \, $\Ker \rho ^L =  P \diag (A ,\, I) K^n_{\infty}(s).$ \,
Hence if
\[
 G =  L_1^{-1} \Theta P \diag (A , \, I)
\]
then \eqref{e.kerincl} is equivalent to
$G \in
s^{-1} K^{n_1 \times n}_{\infty}(s)$.
From \eqref{e.AB} and
\[
 \Sigma = \diag (A, \, 0) + \diag (0, \, B)
\]
we obtain
\[
   L_1^{-1} \Theta L =
                        G \diag (I, \, 0) Q \, + \,
     L_1^{-1} \Theta P \diag (O, I) P^{-1} L.
\]
Now choose
\[
      \Psi =  - G  \diag (I,\, O) Q.
\]
Then $\Psi$ satisfies \eqref{e.Psi} and if we put \,
$ \Theta _1 =  L_1^{-1} \Theta L + \Psi L$ \,
then we have  \linebreak $\Theta _1  \in   K^{n_1 \times n}_{\infty}(s)$,
which proves \eqref{e.L1Psi}. \blackbox

\vspace{2.5mm}

We extend the map $\rho ^{L_1}$ to $ K^n(s)$ and define
\[
       \rho _e ^{L_1} = L_1 \pi _+ L_1^{-1}w, \,\, w \in  K^n(s).
\]

\begin{theorem} \label{Th.main}
The map $\phi : U^L \to U^{L_1}$ is a $K_{\infty}(s)$-module homomorphism
if and only if there exist matrices \,
 $\Theta, \Theta _1 \in K^{n_1 \times n}_{\infty}(s)$ such that
\begin{equation} \label{e.TLwieder}
\Theta L = L_1 \Theta _1
\end{equation}
and
\begin{equation} \label{e.hom}
  \phi \bar{x} = \rho_e ^{L_1} \Theta \bar{x}, \,\,  \bar{x} \in U^L.
\end{equation}
If \eqref{e.TLwieder} holds then we have
\begin{equation} \label{e.ext}
      \rho _e ^{L_1} \Theta \bar{x} =  \rho  ^{L_1} \Theta x
\end{equation}
for all $x \in K^n_{\infty}(s)$.
\end{theorem}

\emph{Proof.}
Let us show first that \eqref{e.TLwieder} implies \eqref{e.ext}. We have
\begin{multline}
    \rho _e ^{L_1} \Theta \bar{x} = L_1 \pi _+ L_1^{-1} \Theta \bar{x} =
  L_1 \pi _+ \Theta _1 L^{-1} \bar{x} = \\
  L_1 \pi _+ \Theta _1 L^{-1} x =  L_1 \pi _+ L_1^{-1} \Theta x =
\rho  ^{L_1} \Theta x.
\end{multline}
Now let  $\phi : U^L \to U^{L_1}$ be a $K_{\infty}(s)$-module homomorphism.
Define       $\Phi = \phi \rho ^L$ such that
\begin{equation} \label{e.pphi}
            \Phi x = \phi \bar{x}, \,\, x \in K^n_{\infty}(s).
\end{equation}
Then \, $     \Phi :  K^n_{\infty}(s) \to U^{L_1}$ \,
is also a $K_{\infty}(s)$-module homomorphism. Because due to Lemma \ref{La.toL1}
there exists a $\tilde{\Theta} \in  K^{n_1 \times n}_{\infty}(s)$ such that
\begin{equation} \label{e.Phix}
\Phi x = \rho  ^{L_1} \tilde{\Theta} x.
\end{equation}
It follows from \eqref{e.pphi} that
 $x, v \in  K^n_{\infty}(s)$ and $\bar{x} = \bar{v}$ imply
$\rho  ^{L_1} \tilde{\Theta} x = \rho  ^{L_1} \tilde{\Theta} v$.
Therefore  we obtain
\begin{equation} \label{e.Thtilde}
\tilde{\Theta} \Ker \rho ^L \subseteq \Ker \rho ^{L_1 }.
\end{equation}
We can replace $\tilde{\Theta}$ in \eqref{e.Phix} and \eqref{e.Thtilde}
by $\Theta = \tilde{\Theta} + L_1 \Psi$
if $\Psi \in s^{-1} K^{n_1 \times n}_{\infty}(s)$ and
 $L_1 \Psi \in  K^{n_1 \times n}_{\infty}(s)$.
From Lemma \ref{La.TT1} we know that starting from \eqref{e.Thtilde}
 we can find a $\Psi$ which yields \eqref{e.TLwieder} with
$\Theta _1 \in  K^{n_1 \times n}_{\infty}(s)$. Thus we have shown
that
\[ \phi \bar{x} = \rho ^{L_1} \Theta x = \rho _e ^{L_1} \Theta \bar{x}
\]
with $\Theta$ satisfying a relation \eqref{e.TLwieder}.

Conversely, if a map $\phi : U^L \to  U^{L_1}$ is defined by
\eqref{e.TLwieder}  and
\eqref{e.hom} then it is easy to verify that  $\phi$ is a
$K_{\infty}(s)$-module homomorphism. \blackbox

\vspace{2.5mm}
We remark that Theorem \ref{Th.main} remains true if condition
\eqref{e.TLwieder} is replaced by
\[
    \pi _+ L_1 ^{-1} \Theta =  \pi _+ \Theta _1 L^{-1}.
\]

Given the duality \eqref{e.dual}
between $U^L$ and $U^{L^T}$ it is not difficult to obtain the dual map of
$\phi$.
We set $\bar{\bar{w}} = \rho ^{L_1^T}w, \, w \in
K^{n_1}_{\infty}(s)$.

\begin{theorem} \label{Th.dual}
Let  \,  $\Theta, \Theta _1 \in K^{n_1 \times n}_{\infty}(s)$  \,
be such that  \,\,$\Theta L \, = \, L_1 \Theta _1 $. \,
Let \linebreak  $\phi :  U^L \to U^{L_1}$ be defined by \eqref{e.hom}. Then
the dual map
\[
    \phi^* : U^{L_1^T} \to  U^{L^T}
\]
is given by
\[
 \phi^* \bar{\bar{w}} = \rho ^{L^T} \Theta _1^T w, \,\,
\bar{\bar{w}} \in  U^{L_1^T}.
\]
\end{theorem}

We now turn to surjectivity and injectivity.
For a pair  $\Theta \in  K^{n_1 \times n}_{\infty}(s)$ and
$L_1 \in K^{n_1 \times n_1}$ we set
$  (\Theta, \, s^{-1}L_1)_l = I $
if there exist proper rational matrices  $C$ and $D$ such that
\begin{equation}  \label{e.left}
                     \Theta C +  s^{-1}L_1 D = I.
\end {equation}
Similarly, for $\Theta _1 \in   K^{n_1 \times n}_{\infty}(s)$ and
$L \in K^{n \times n}$ we write
$       (\Theta _1, s^{-1}L)_r =I $
if  $(\Theta _1 ^T, s^{-1}L^T)_l = I$.

\begin{theorem} \label{Thm.surj}
Let $\phi : U^L \to U^{L_1}$ be defined by \eqref{e.ext}
 and
\eqref{e.TLwieder}. Then
\begin{enumerate}
\item[\rm{(i)}]
 $\phi$ is surjective if and only if \,
      $ (\Theta, \, s^{-1}L_1)_l = I $,
\item[\rm{(ii)}]
$\phi$ is injective if and only if \,
     $ (\Theta _1, s^{-1}L)_r =I$.
\end{enumerate}
\end{theorem}

\emph{Proof.} (i) Assume first that  $\phi$ is surjective. Let $w \in
 K^{n_1}_{\infty}(s)$ be given. Then $\rho ^{L_1}w =
\rho ^{L_1} \Theta v$ for some $v \in  K^{n}_{\infty}(s)$. We have
$w - \Theta v \in \Ker \rho ^{L_1}$, which implies
   \[
       w \in \Theta  K^{n}_{\infty}(s) +  s^{-1}L_1  K^{n}_{\infty}(s)
\]
or equivalently $ (\Theta, \, s^{-1}L_1)_l = I $.  Conversely, suppose
that \eqref{e.left} holds. To show that
$w = \rho ^{L_1}x$ is in $\phi U^L$ we note that
\eqref{e.left} implies $x = \Theta v +  s^{-1}L_1 x_2$ for some
$v \in K^{n}_{\infty}(s), \, x_2 \in  K^{n_1}_{\infty}(s)$. Because of
$ s^{-1}L_1 x_2 \in \Ker \rho ^{L_1}$ we obtain $w = \rho ^{L_1} \Theta v =
\phi \bar{v}$.\\
(ii) By duality the statement follows at once from (i).
\blackbox

If $M$ is a finitely generated $p$-module over a principal ideal
domain and $S$ is a submodule and $Q$ is a quotient module of $M$ then
the relations between the invariants of $M$ and those of $S$ and $Q$
are well known (see e.g. \cite[p. 92, 93]{Ri}). We complete our
note with a corresponding  observation on the existence of surjective and
injective homomorphisms. Let
\[ s^{- \alpha _1}, \dots ,  s^{- \alpha _t}, \,\,
\alpha _1 \geq \dots \geq \alpha _t,
\]
and
 \[s^{- \gamma _1}, \dots ,  s^{- \gamma _p}, \,\,
\gamma _1 \geq \dots \geq \gamma _p,
\]
be the infinite elementary divisors of $L$ and  $L_1$, respectively.
Then there exists a surjective $ K_{\infty}(s)$-module
homomorphism $\phi : U^L \to U^{L_1}$ if and only if
\[
    t \geq p \quad  \mbox{and} \quad  \alpha _1 \geq \gamma _1,  \dots
,
\alpha _p \geq \gamma _p,
\]
and there exists an injective $\phi$
if and only if
\[
  t \leq p \quad \mbox{and} \quad  \alpha _1 \leq \gamma _1,  \dots ,
\alpha _t \leq \gamma _t.
\]

\end{document}